\documentclass[12pt,a4paper]{amsart}
\usepackage{amscd,amssymb,amsopn,amsmath,amsthm,graphics,amsfonts,enumerate,verbatim,calc}
\headsep=1cm \numberwithin{equation}{section}
\numberwithin{equation}{section}

\begin{document}
\pagenumbering{gobble}

\begin{flushright}
\end{flushright}

\title{On Poisson Distribution}
\author{$^{*}$Dharmendra Kumar Singh}
\address{Department of Mathematics, University Institute of Engineering and Technology, CSJM University, Kanpur, India}
\email{drdksinghabp@gmail.com}


\begin{abstract}
The object of this paper is to study and develop a Poisson distribution in generalized Wright function form.\vspace{0.2 cm} 
\end{abstract}

\maketitle
\section{Introduction}
{\textbf{Generalized Wright Function}}\\

The generalized Wright hypergeometric function $_{p}\Psi_{q}$ [2, 4] is defined by the series
\begin{equation}
_{p}\Psi_{q}(z)=_{p}\Psi_{q}\left[z\vert \begin{array}{c} \left(a_{i}, \alpha_{i}\right)_{1, p}\\ \left(b_{j}, \beta_{j}\right)_{1, q}\end{array}\right]=\sum_{k=0}^{\infty}\frac{\prod_{i=1}^{p} \Gamma(a_{i}+\alpha_{i}k)}{\prod_{j=1}^{q}\Gamma(b_{j}+\beta_{j}k)}\frac{z^{k}}{k!}.
\end{equation}
Here $a_{i}, b_{j}\in C$ and $\alpha_{i}, \beta_{j}\in R \left(i= 1, 2, ..., p; j= 1, 2, ..., q\right)$. Conditions for the existence of the generalized Wright in terms of $H-$ function were established in [1]. The particular $_{p}\Psi_{q}(z)$ is an entire function if there holds the condition
\begin{equation}
\sum_{j=1}^{q}\beta_{j}-\sum_{i=1}^{p} \alpha_{i}>1.
\end{equation}

{\textbf{Mittag Leffler Function}}\\
 The Mittag-Leffler function $E_{\alpha}(z)$ [5] is defined by
 \begin{equation}
 E_{\alpha}(z)=\sum_{k=o}^{\infty}\frac{z^{k}}{\Gamma (1+\alpha k)}.
 \end{equation} 
 For $\alpha=1$,  (1.3) arrive at exponential function $e^{x}$. In 1905 Wiman's [3] arose generalization of the Mittag-Leffler function (1.3)
 \begin{equation}
 E_{\alpha, \beta}(z)=\sum_{k=o}^{\infty}\frac{z^{k}}{\Gamma (\alpha k+\beta)}.
 \end{equation}
 Here $E_{\alpha, 1}(z)=E_{\alpha}(z)$.
Prabhakar [7] introduced the  Mittag-Leffler function $E_{\alpha, \beta}^{\gamma}(z)$ defined by
\begin{equation}
E_{\alpha, \beta}^{\gamma}(z)=\sum_{k=o}^{\infty}\frac{(\gamma)_{k} }{\Gamma (\alpha k+\beta)} \frac{z^{k}}{k!},
\end{equation}
where $\alpha,\beta, \gamma \in C, \Re(\alpha)>0, \Re(\beta)>0, \Re(\gamma)>0$ and $(\gamma)_{n}$ is the Pochhamer symbol defined by\\
\begin{equation}
(\gamma)_{n}:=\left\{\begin{array}{ccc} 1~~~~~~~~~~~ &(n=0)\\ \gamma (\gamma+1) (\gamma+2)...(\gamma+n-1)~~~~~~~~~~~&(n\in N)\\ \frac{\Gamma(\gamma+n)}{\Gamma(\gamma)}~~~~~~~~~~&(\gamma\in C)\end{array}\right\}
\end{equation}
{\textbf{Poisson Distribution}}\\
The series
\[1+m+\frac{m^{2}}{2!}+\frac{m^{3}}{3!}+... =\sum_{r=0}^{\infty}\frac{m^{r}}{r!}.\]
Converges, for all values of $m$, to $e^{m}$. Consider the function $f(r)$ defined  by
\[f(r)=\frac{m^{r}e^{-m}}{r!}, ~~~~~~r= 0, 1, 2,...,\]
\begin{equation}
= 0 \bigskip~~~~~~~~~~~~~~~~{\text{elsewhere}},
\end{equation}
where $m>0$. Since $m>0$, then $f(r)\geq 0$ and\\
\[\sum_{r}f(r)=\sum_{r=0}^{\infty}\frac{m^{r}e^{-m}}{r!}=e^{-m}\sum_{r=0}^{\infty}\frac{m^{r}}{r!}\]
\begin{equation}
 =e^{-m}e^{m}=1
\end{equation}
that is, $f(r)$ satisfies the conditions of being a probability density function of a discrete type of random variable. A random variable $X$ which has a probability density function in the form of $f(r)$ is said to have a Poisson distribution, and any such $f(r)$ is called a Poisson probability distribution function [6].\\
{\textbf{Wright-type Poisson Distribution}}\\
The topic which I am introducing is the outline for this research paper
\begin{equation}
P(\alpha, \beta, 1, m, r)=\frac{m^{r}}{{_{1}\Psi_{1}}\left[m \vert \begin{array}{c}(1, 1)\\(\beta, \alpha)\end{array}\right]\Gamma(\alpha r+\beta )}.
\end{equation}
\begin{equation}
P(\alpha, \beta, 1, m, r)=\frac{m^{r}}{E^{1}_{\alpha, \beta}(m)\Gamma(\alpha r+\beta)}.
\end{equation}
Where $\alpha>0, \beta>0$ and $m>0$. Note that for ${\alpha=\beta=1}$ we arrive at (1.7).
\section{Main Results}
{\textbf{Probability Density Function}}\\
Since 
\begin{equation*}
\sum_{r=0}^{\infty}P(\alpha, \beta, 1, m, r)=\sum_{r=0}^{\infty}\frac{m^{r}}{{_{1}\Psi_{1}}\left[m \vert \begin{array}{c}(1, 1)\\(\beta, \alpha)\end{array}\right]\Gamma(\alpha r+\beta)}.
\end{equation*}
\begin{equation*}
\sum_{r=0}^{\infty}P(\alpha, \beta, 1, m, r)=\sum_{r=0}^{\infty}\frac{m^{r}}{\sum_{r=0}^{\infty}\frac{\Gamma(1+r) m^{r}}{\Gamma(\alpha r+\beta)r!} \Gamma(\alpha r+\beta)}=1.
\end{equation*}
Therefore (1.9) and (1.10) is a probability density functions.\\
{\vspace{4 mm}}\\
{\text{\large{\bf{Mean:}}}\\
{\vspace{1 mm}}\\
{\textbf{Method I}}
\begin{equation*}
\sum_{r=0}^{\infty} r P(\alpha, \beta, m, r)=\sum_{r=0}^{\infty}\frac{r m^{r}}{{_{1}\Psi_{1}}\left[m \vert \begin{array}{c}(1, 1)\\(\beta, \alpha)\end{array}\right]\Gamma(\alpha r+\beta)}.
\end{equation*}
\begin{equation*}
=\sum_{r=0}^{\infty}\frac{(r+1-1)\Gamma(r+1) m^{r}}{{_{1}\Psi_{1}}\left[m \vert \begin{array}{c}(1, 1)\\(\beta, \alpha)\end{array}\right]\Gamma(\alpha r+\beta )r!}.
\end{equation*}
\begin{equation*}
=\frac{1}{{_{1}\Psi_{1}}\left[m \vert \begin{array}{c}(1, 1)\\(\beta, \alpha)\end{array}\right]}\left[\sum_{r=0}^{\infty}\frac{(r+1)\Gamma(r+1)m^{r}}{\Gamma(\alpha r+\beta)r!}-\sum_{r=0}^{\infty}\frac{\Gamma(r+1) m^{r}}{\Gamma(\alpha r+\beta)r!}\right]
\end{equation*}
\begin{equation*}
=\frac{1}{{_{1}\Psi_{1}}\left[m \vert \begin{array}{c}(1, 1)\\(\beta, \alpha)\end{array}\right]}\left[\sum_{r=0}^{\infty}\frac{\Gamma(r+2)m^{r}}{\Gamma(\alpha r+\beta)r!}-\sum_{r=0}^{\infty}\frac{\Gamma(r+1) m^{r}}{\Gamma(\alpha r+\beta)r!}\right]
\end{equation*}
\begin{equation}
=\frac{1}{{_{1}\Psi_{1}}\left[m\vert\begin{array}{c}(1, 1)\\ (\beta, \alpha)\end{array}\right]}\left[{_{1}\Psi_{1}}\left[m\vert\begin{array}{c}(2, 1)\\ (\beta, \alpha)\end{array}\right]-{_{1}\Psi_{1}}\left[m\vert\begin{array}{c}(1, 1)\\(\beta, \alpha)\end{array}\right]\right]
\end{equation}
{\vspace{5 mm}}\\
{\bf{Method II}}
\[\sum_{r=0}^{\infty} r P(\alpha, \beta, m, r)=\frac{1}{\alpha}\sum_{r=0}^{\infty}\frac{\alpha r m^{r}}{{_{1}\Psi_{1}}\left[m \vert \begin{array}{c}(1, 1)\\(\beta, \alpha)\end{array}\right]\Gamma(\alpha r+\beta)}.\]
\[=\frac{1}{\alpha}\sum_{r=0}^{\infty}\frac{(\alpha r+\beta-1+1-\beta) m^{r}(1)_{r}}{{_{1}\Psi_{1}}\left[m \vert \begin{array}{c}(1, 1)\\(\beta, \alpha)\end{array}\right]\Gamma(\alpha r+\beta)~~~r!}.\]
\[=\frac{1}{\alpha}\frac{1}{{_{1}\Psi_{1}}\left[m \vert \begin{array}{c}(1, 1)\\(\beta, \alpha)\end{array}\right]}\left[\sum_{r=0}^{\infty}\frac{m^{r}(1)_{r}}{\Gamma(\alpha r+\beta-1)~~~r!}+(1-\beta)\sum_{r=0}^{\infty}\frac{m^{r}(1)_{r}}{\Gamma(\alpha r+\beta)~~~r!}\right]\]
\begin{equation}
=\frac{1}{\alpha}\frac{1}{{_{1}\Psi_{1}}\left[m\vert\begin{array}{c}(1, 1)\\ (\beta, \alpha)\end{array}\right]}\left[ E^{1}_{\alpha, \beta-1}(m)+(1-\beta)E^{1}_{\alpha, \beta}(m)\right]
\end{equation}
\begin{equation}
=\frac{1}{\alpha}\frac{1}{{_{1}\Psi_{1}}\left[m\vert\begin{array}{c}(1, 1)\\ (\beta, \alpha)\end{array}\right]}\left[ {_{1}\Psi_{1}}\left[m\vert\begin{array}{c}(1, 1)\\ (\beta-1, \alpha)\end{array}\right]+(1-\beta){_{1}\Psi_{1}}\left[m\vert\begin{array}{c}(1, 1)\\ (\beta, \alpha)\end{array}\right]\right]
\end{equation}
Note that for $\alpha=\beta=1$ above expressions (2.1), (2.2) and (2.3) arrive mean of Poisson distribution [6].\\
{\vspace{1 mm}}\\
\textbf{Variance:}
{\vspace{1 mm}}\\
{\bf{Method I}}
\[\sum_{r=0}^{\infty}r^{2}P(\alpha, \beta, 1,  m, r)=\sum_{r=0}^{\infty}\frac{r^{2} m^{r}}{{_{1}\Psi_{1}}\left[m \vert \begin{array}{c}(1, 1)\\ (\beta, \alpha)\end{array}\right]\Gamma(\alpha r+\beta)}.\]
\[=\sum_{r=0}^{\infty}\frac{(r(r-1)+r)(1)_{r} m^{r}}{{_{1}\Psi_{1}}\left[m \vert \begin{array}{c}(1, 1)\\(\beta, \alpha)\end{array}\right]\Gamma(\alpha r+\beta )~~~r!}.\]
\[=\frac{1}{{_{1}\Psi_{1}}\left[m \vert \begin{array}{c}(1,1)\\(\beta, \alpha)\end{array}\right]}\left[\sum_{r=2}^{\infty}\frac{r(r-1)\Gamma(r-1)}{\Gamma(r-1)\Gamma(\alpha r+\beta)}+\sum_{r=0}^{\infty}\frac{r }{\Gamma(\alpha r+\beta)}\right]\frac{(1)_{r} m^{r}}{r!}\]
\[=\frac{1}{{_{1}\Psi_{1}}\left[m \vert \begin{array}{c}(1,1)\\(\beta, \alpha)\end{array}\right]}\left[\sum_{r=2}^{\infty}\frac{\Gamma(r+1) \Gamma(r+1) m^{r}}{\Gamma(r-1)\Gamma(\alpha r+\beta) ~~~r!}+\sum_{r=0}^{\infty}\frac{(r+1-1)\Gamma(r+1) m^{r}}{\Gamma(\alpha r+\beta)~~~r!}\right]\]
\[=\frac{1}{{_{1}\Psi_{1}}\left[m \vert \begin{array}{c}(1, 1)\\(\beta, \alpha)\end{array}\right]}\left[\sum_{r=2}^{\infty}\frac{\Gamma(r+1)\Gamma(r+1)m^{r}}{\Gamma(r-1)\Gamma(\alpha r+\beta)r!}+\sum_{r=0}^{\infty}\frac{\Gamma(r+2) m^{r}}{\Gamma(\alpha r+\beta)r!}-\sum_{r=0}^{\infty}\frac{\Gamma(r+1) m^{r}}{\Gamma(\alpha r+\beta)r!}\right]\]
\begin{equation}
=\frac{1}{\frac{1}{_{1}\Psi_{1}}\left[m\vert\begin{array}{c}(1, 1)\\ (\beta, \alpha)\end{array}\right]}\left[{_{2}\Psi_{2}}\left[m\vert\begin{array}{cc}(1, 1) (1, 1)\\(-1, 1) (\beta, \alpha)\end{array}\right]+{_{1}\Psi_{1}}\left[m\vert\begin{array}{c}(2, 1)\\(\beta, \alpha)\end{array}\right]-
{_{1}\Psi_{1}}\left[m\vert\begin{array}{c}(1, 1)\\(\beta, \alpha)\end{array}\right]\right].
\end{equation}
Note that for $\alpha=\beta=1$ we arrive variance of Poisson distribution [6].\\
{\vspace{1 mm}}\\
{\bf{Method II}}
\[\sum_{r=0}^{\infty}r^{2}P(\alpha, \beta, 1,  m, r)=\sum_{r=0}^{\infty}\frac{r^{2} m^{r}}{{_{1}\Psi_{1}}\left[m \vert \begin{array}{c}(1, 1)\\ (\beta, \alpha)\end{array}\right]\Gamma(\alpha r+\beta)}.\]
\[=\frac{1}{\alpha^{2}}\sum_{r=0}^{\infty}\frac{\alpha^{2} r^{2} (1)_{r} m^{r}}{{_{1}\Psi_{1}}\left[m \vert \begin{array}{c}(1, 1)\\ (\beta, \alpha)\end{array}\right]\Gamma(\alpha r+\beta)~~~r!}.\]
\[=\frac{1}{\alpha^{2}{_{1}\Psi_{1}}\left[m \vert \begin{array}{c}(1, 1)\\ (\beta, \alpha)\end{array}\right]}\sum_{r=0}^{\infty}\frac{\left[(\alpha r+\beta-1)(\alpha r+\beta-2)+(\alpha r+\beta-1)(3-2\beta)+(1-\beta)^{2}\right] (1)_{r} m^{r}}{\Gamma(\alpha r+\beta)~~~r!}.\]
\[=\frac{1}{\alpha^{2}{_{1}\Psi_{1}}\left[m \vert \begin{array}{c}(1, 1)\\ (\beta, \alpha)\end{array}\right]}\left[\sum_{r=0}^{\infty}\frac{(1)_{r} m^{r}}{{\Gamma(\alpha r+\beta-2)~~~r!}}+\sum_{r=0}^{\infty} \frac{(3-2\beta) (1)_{r} m^{r}}{{\Gamma(\alpha r+\beta-1)~~~r!}}+\sum_{r=0}^{\infty}\frac{(1-\beta)^{2} (1)_{r} m^{r}}{\Gamma(\alpha r+\beta)~~~r!}\right].\]
\begin{equation}
=\frac{1}{\alpha^{2}{_{1}\Psi_{1}}\left[m \vert \begin{array}{c}(1, 1)\\ (\beta, \alpha)\end{array}\right]}\left[E^{1}_{\alpha, \beta-2}{m}+(3-2\beta) E^{1}_{\alpha, \beta-1}(m)+(1-\beta)^{2}E^{1}_{\alpha, \beta}\right]
\end{equation}
\begin{equation*}
=\frac{1}{\alpha^{2}{_{1}\Psi_{1}}\left[m \vert \begin{array}{c}(1, 1)\\ (\beta, \alpha)\end{array}\right]}{_{1}\Psi_{1}}\left[m \vert \begin{array}{c}(1, 1)\\ (\beta-2, \alpha)\end{array}\right]+(3-2\beta) {_{1}\Psi_{1}}\left[m \vert \begin{array}{c}(1, 1)\\ (\beta-1, \alpha)\end{array}\right]
\end{equation*}
\begin{equation}
\left. +(1-\beta)^{2}{_{1}\Psi_{1}}\left[m \vert \begin{array}{c}(1, 1)\\ (\beta, \alpha)\end{array}\right]\right\}
\end{equation}
Note that for $\alpha=\beta=1$, (2.4), (2.5) and (2.6) arrive variance of Poisson distribution [6].\\
{\vspace{2 mm}}\\
\textbf{Recurrence Relation}\\
{\vspace{0.5 mm}}\\
Since (1.9)
\begin{equation}
P(\alpha, \beta, 1, m, r)=\frac{ m^{r}}{{_{1}\Psi_{1}}\left[m \vert \begin{array}{c}(1, 1)\\ (\beta, \alpha)\end{array}\right]\Gamma(\alpha r+\beta )}.
\end{equation}
 and
\begin{equation}
P(\alpha, \beta, 1, m, r+1)=\frac{ m^{r+1}}{{_{1}\Psi_{1}}\left[m \vert \begin{array}{c}(1, 1)\\ (\beta, \alpha)\end{array}\right]\Gamma(\alpha r+\alpha+\beta)}.
\end{equation}
From (2.7)and (2.8)
\begin{equation}
P(\alpha, \beta, 1, m, r+1)=\frac{m\Gamma(\alpha r+\beta)}{\Gamma(\alpha r+\alpha+\beta)} P(\alpha, \beta, 1, m, r).
\end{equation}

Note that for $\alpha=\beta=1 $ above expression (2.9) arrive at reccurence relation of Poisson distribution [6].\\
{\vspace{2 mm}}\\
\textbf{Moment generating Function}\\
{\vspace{0.5 mm}}
Let $X$ be a Poisson variate Then MGF of $X$ is defined by
\begin{equation}
M_{x}(t)=E(e^{tx}).
\end{equation}
\[M_{x}(t)=\frac{1}{\sum_{n=0}^{\infty}\frac{m^{n}\Gamma(1+ n)}{\Gamma(\alpha n+\beta) n!}}\sum_{r=0}^{\infty}\frac{e^{tr}m^{r}}{\Gamma(\alpha r+\beta)}.\]
\[=\frac{E_{\alpha, \beta}(e^{t}m)}{E_{\alpha, \beta}^{1}(m)}.\]
\begin{equation}
=\frac{E_{\alpha, \beta}(e^{t}m)}{{_{1}\Psi_{1}}\left[m\vert\begin{array}{c}(1, 1)\\ (\beta, \alpha)\end{array}\right]}.
\end{equation}
Note that for $\alpha=\beta=1$ above expression (2.11) arrive at Moment generating function of Poisson distribution [6].\\



\


\end{document}